\documentclass[preprint,12pt]{elsarticle}

\usepackage{amssymb}
\usepackage{amsmath}

\begin{document}

\begin{frontmatter}



\title{Exact solution to a class of problems for the Burgers' equation on bounded intervals}

 \author[label1]{Kwassi Anani}
 \affiliation[label1]{organization={Department of Mathematics, University of Lomé},
             addressline={kanani@univ-lome.tg},
             city={Lomé},
             postcode={01 BP 1515},
             country={Togo}}

\author[label2]{Mensah Folly-Gbetoula}
 \affiliation[label2]{organization={School of Mathematics, University of the Witwatersrand},
             addressline={Mensah.Folly-Gbetoula@wits.ac.za},
             city={Johannesburg},
             postcode={Wits 2050},
             country={South Africa}}

%

\begin{abstract}
In this study, we consider Burgers' equation with fixed Dirichlet boundary conditions on generic bounded intervals. By employing the Hopf-Cole  transformation and a recently established exact operational solution for linear reaction-diffusion equations, an exact solution in the time domain is derived through inverse Laplace transforms. In the event that analytic inverses do in fact exist, they can be obtained in closed form through the use of Mellin transforms. Nevertheless, highly efficient algorithms are available, and numerical inverses in the time domain are always feasible, regardless of the complexity of the Laplace domain expressions. Two illustrative tests demonstrate that the results align closely with those of classical exact solutions. In comparison to the solutions obtained with series expressions or by numerical methods, closed-form expressions, even in the Laplace domain, represent a novel alternative, offering new insights and perspectives. The exact solution via the inverse Laplace transform is shown to be more computationally efficient, providing a reference point for numerical and semi-analytical methods.
\end{abstract}


\begin{highlights}
\item Burgers' equation is considered on generic bounded intervals
\item Exact solution in the time domain is derived by means of inverse Laplace transforms
\item The closed-form expression of the solution offers new alternative and perspectives
\item The solution provides a reference point for numerical and semi-analytical methods
\end{highlights}

\begin{keyword}
Burgers' equation \sep Hopf-Cole  transformation \sep integral transforms \sep closed-form solution \sep Laplace domain \sep Numerical efficiency


\end{keyword}

\end{frontmatter}




\section{Introduction}\label{sec:1}

The numerical simulation of linear and nonlinear convective and diffusive partial differential equations is of significant importance, particularly in the context of the Burgers' equation (see e.g. Ba\c{s}han 2022, \cite{Ba-nonlinear} and Kumar et al. \cite{Ku-Ku-evolutionary}). This equation is used in fluid dynamics to model viscous fluids, study shock waves, turbulence, and other flow characteristics, and examine shocks in inviscid flows (see \cite{Kh-Iq-mathematical}, \cite{Ha-Gu-anomalies1} and \cite{Ha-Gu-anomalies2} for examples). Furthermore, Burgers' equation is used to model traffic flow, describing the evolution of vehicular density. It optimizes traffic patterns, improves the understanding of congestion, and can permit to design superior road systems (see e.g.  \cite{Zh-Wu-continuous} and \cite{Zh-al-cooperative}). The equation is also used in image and signal processing, particularly for noise reduction. It allows to reduce noise while maintaining edges in signals or images. This technique is crucial in domains like medical imaging, where noise reduction is key for accurate diagnosis (see e.g. \cite{Gu-Ji-modeling} and \cite{Ch-al-sonification}). The behavior of non-linear waves in a variety of media, including acoustics, water waves, and plasma, is frequently modeled using the Burgers' equation. For example, it can be used to elucidate shock waves in gases and other fluids (see references \cite{Ru-He-quadratically} and \cite{Ku-al-dynamics}). In acoustics, it facilitates comprehension of sound wave propagation in dissipative media. For a review of the applications of Burgers' equation in mathematical physics, see Bonkile et al. \cite{Bo-Aw-systematic}.

The first historical solutions in terms of the exact ratio of infinite series were obtained through the Hopf-Cole  transformation, as in \cite{Ho-partial, Co-quasilinear}. A review of exact solutions to the one-dimensional Burgers' equation with diverse boundary and initial conditions has been conducted by Benton and Platzman in \cite{Be-Pl-table}. In recent decades, numerical methods have become the predominant approach for studying Burgers' equations. A variety of numerical methods may be cited, including finite difference, finite element, spectral element, spline approximations, modified cubic B-splines collocation, hybrid numerical scheme, and so forth. Further examples are referenced in the work of Dhawan et al. \cite{Dh-Ka-contemporary}. Also, semi-analytical methods have been utilized in the context of Adomian decomposition, as demonstrated by Naghipour and Manafian \cite{Na-Ma-mathematical} and in the approach of Zeidan et al. \cite{Ze-Ch-mathematical}. Homotopy perturbation techniques have been employed by Lal and Yadav \cite{La-Ya-approximate}, while variational iterations have been utilized by Biazar and Aminikhah \cite{Bi-Am-exact}. Additionnally, the Lie symmetry technique has also been utilized for the symmetry reductions of coupled Burgers equations as in \cite{Os-al-approximate}.

From one perspective, numerical methods are undergoing continual improvement in order to enhance their capacity to adapt to fluctuations in problem parameters. Conversely, semi-analytical methods do not always yield a closed-form solution, which can result in the consumption of considerable resources and CPU time. The objective of this paper is to derive a closed-form operational solution and subsequently deduce the exact solution in the time domain through the application of the inverse Laplace transform operator. In instances where exact analytic inverses into the time domain are unavailable, computationally efficient algorithms exist that facilitate the numerical implementation of inverse Laplace transforms. The paper is organized as follows. Section \ref{sec:2} introduces the problem and presents a solution utilizing the Hopf-Cole  transformation and a previously derived operational solution in a closed-form expression for the reaction-diffusion equation on finite intervals. A comparison of the present solutions with classical exact and numerical results is presented in Section \ref{sec:3}. The paper concludes with a summary and an outline of future research directions.
\section{Exact implicit solution}\label{sec:2}
Following examples as in \cite{Gu-al-fifth, Se-accurate}, the viscid Burgers' equation is now considered on a generic bounded interval of the real line and consists to determine the function $ w(x,t) $ as:
\begin{equation}
	w_{t}-a^2 w_{xx}+ ww_{x}=0 \quad l_1 < x < l_2,\quad 0 < t \le T,
	\label{eq01}
\end{equation}
subject to a non-homogeneous initial condition:
\begin{equation}
	w(x,0)=w_{0}(x), \quad    l_1 \le x \le l_2,
	\label{eq02}
\end{equation}
and to the following Dirichlet boundary conditions for  $0 \le t \le T $,
\begin{equation}
	w(l_1,t)=\alpha_1,
	\label{eq03}
\end{equation}
\begin{equation}
	w(l_2,t)=\alpha_2,
	\label{eq04}
\end{equation}
where $\alpha_1$ and $\alpha_2$ are two arbitrary real numbers and $t > 0 $ when $ T= \infty$. The coefficient of kinematic viscosity $a^2$ is defined by $a^2=1/Re$, where $Re$ is the Reynolds number. The Burgers' equation can be considered as a simplified form of the Navier-Stokes equations. This equation is a balance between time evolution, nonlinearity, and diffusion due to the simultaneous presence of a nonlinear convective term ($ ww_{x}$) and a diffusive term ($a^2 w_{xx}$). In this sense, Burgers' equation is a suitable benchmark problem for the development and testing of numerical methods, due to its equilibrium between nonlinear convection and diffusion terms. It should be noted that the equation is parabolic when $a^2 \neq 0$ and hyperbolic when $a^2 = 0$. As $a^2$ approaches zero, equation (\ref{eq01}) becomes inviscid Burgers' equation, which serves as a model for nonlinear wave propagation.

The classical Hopf-Cole  transformation is applied under the condition that the new function, designated as $u(x,t)$, does not exhibit cancellation over the domain $[l_1, l_2] \times [0, T]$. We pose
\begin{equation}
	w(x,t)=-2a^2\frac{u_{x}(x,t)}{u(x,t)},
	\label{eq05}
\end{equation}
and the system (\ref{eq01})-(\ref{eq04}) is transformed unto a boundary value problem for the one-dimensional heat diffusion equation:
\begin{equation}
	u_{t}(x,t)-a^2 u_{xx}(x,t)=0, \quad l_1 < x < l_2,\quad 0 < t \le T,
	\label{eq06}
\end{equation}
subject to the initial condition
\begin{equation}
	\displaystyle u(x,0)=u(l_1,0)\exp\left( -\frac{1}{2a^2}\int_{l_1}^{x}w_{0}(y)dy\right)=\varphi(x),  \quad    l_1 \le x \le l_2,
	\label{eq07}
\end{equation}
and to the homogeneous Robin boundary conditions for $0 \le t \le T $
\begin{equation}
	\alpha_1 u(l_1,t)+2a^2 u_{x}(l_1,t)=0,
	\label{eq08}
\end{equation}
\begin{equation}
	\alpha_2 u(l_2,t)+2a^2 u_{x}(l_2,t)=0.
	\label{eq09}
\end{equation}
Equation (\ref{eq06}) is subject to the initial condition (\ref{eq07}) (it is assumed that $u(l_1,0)$ is not null), and to homogeneous Robin boundary conditions (\ref{eq08}) and (\ref{eq09}). Now, the system (\ref{eq06})-(\ref{eq09}) is well defined with a boundary value problem for a linear reaction-diffusion equation with constant coefficients. The exact solution has been provided in the Laplace domain for such problems in Anani \cite{An-approximations}.

We briefly recall that, for a nonhomogeneous linear reaction-diffusion equation with constant coefficients:
\begin{equation}
	\frac{\partial u}{\partial t}-a^2\frac{\partial^{2}u}{\partial x^{2}} + bu =f(x,t), \quad l_1 < x < l_2,\quad 0 < t \le T,
	\label{eq10}
\end{equation}
subject to the initial condition
\begin{equation}
	u(x,0)=\varphi(x), \quad    l_1 \le x \le l_2,
	\label{eq11}
\end{equation}
and to the boundary conditions for $ 0 \le t \le T $:
\begin{equation}
	\alpha_1 u(l_1,t)+\beta_1 \frac{\partial u}{\partial x}(l_1,t)=g_1(t), \quad \alpha_1^2+\beta_1^2 \ne 0,
	\label{eq12}
\end{equation}
\begin{equation}
	\alpha_2 u(l_2,t)+\beta_2 \frac{\partial u}{\partial x}(l_2,t)=g_2(t),  \quad \alpha_2^2+\beta_2^2 \ne 0,
	\label{eq13}
\end{equation}
its analog in the Laplace domain can be written in the form of a two-point boundary value problem as follows:
\begin{equation}
	-a^2\frac{d^2 U}{dx^{2}}(x,p)+(b+p)U(x,p)=F(x,p)+\varphi(x),
	\label{eq14}
\end{equation}
with the boundary conditions:
\begin{equation}
	\alpha_1 U(l_1,p)+\beta_1 \frac{d U}{d x}(l_1,p)=G_1(p), \quad \alpha_1^2+\beta_1^2 \ne 0,
	\label{eq15}
\end{equation}
\begin{equation}
	\alpha_2 U(l_2,p)+\beta_2 \frac{d U}{d x}(l_2,p)=G_2(p),   \quad \alpha_2^2+\beta_2^2 \ne 0.
	\label{eq16}
\end{equation}
In equations (\ref{eq14})-(\ref{eq16}), $p$ is the Laplace domain variable for the time $t$, $U(x,p)$ denotes the Laplace transform of  $ u(x,t) $, and $ U_x $ stands for the derivative $ \frac{dU}{dx} $ of $U(x,p)$. In equations (\ref{eq10})-(\ref{eq13}), we adopt the same approach as in \cite{An-approximations} and assume that the source term, represented by $f(x,t)$, is continuous on the interval $(l_1, l_2)\times (0, T]$, and the function $ \varphi(x) $ satisfies the Dirichlet conditions on the interval $[l_1, l_2]$. That is, this function is piecewise continuous or can be expressed as a convergent series of eigenfunctions, which constitutes a complete basis for the related Sturm-Liouville problem. Additionally, all time-dependent functions are presumed to be of exponential order in relation to $t$, as well as their temporal derivatives where they are involved, particularly when $T = \infty$. For a detailed examination of the requisite conditions for a function to have a Laplace transform, please refer to the following sources: Herron and Foster \cite{Her-partial} or Debnath and Bhatta \cite{De-integral}. Thus, $F(x,p)$, $G_1(p)$ and $G_2(p)$  are respectively the Laplace transform of the source term $ f(x,t) $, and the boundary-related functions $ g_1(t) $ and $ g_2(t) $. Under these hypotheses, the  uniqueness of the solution $u$ has been proven. And, the exact operational solution of the system (\ref{eq10})-(\ref{eq13}), or equivalently, the exact solution of the analog problem (\ref{eq14})-(\ref{eq16}) is:
\begin{equation}
	\begin{array}{ll}
		\displaystyle	U(x,p)= \frac {1}{2\, \sqrt{b+p}} \,{\exp\left( {-{\frac {\left(l_{{2}}-x \right)  \sqrt{b+p}}{a}}}\right)} \left[ { U(l_2,p)} \sqrt{b+p}+a{ U_x(l_2,p)}\right]\\
		\displaystyle+ {\frac {1}{2\, \sqrt{b+p}} {\exp\left( {-{\frac { \left( x-l_{{1}} \right)  \sqrt{b+p}}{a}}}\right) }\left[{ U(l_1,p)} \sqrt{b+p}-a{U_x(l_1,p)}\right]}+R \left( x,p \right).
		\label{eq17}
	\end{array}
\end{equation}
In equation (\ref{eq17}), the four $p$-domain coefficients  $U(l_1,p)$, $ U_x(l_1,p)$, $ U(l_2,p)$ and $ U_x(l_2,p)$ are determined by the following linear system of four equations:
$$
\begin{array}{ll}
	\left[
	\begin{array}{cccc}
		\displaystyle U(l_1,p)\\
		\displaystyle U_x(l_1,p)\\
		\displaystyle U(l_2,p) \\
		\displaystyle U_x(l_2,p)
	\end{array}
	\right]	=\left[
	\begin{array}{cccc}
		\displaystyle \alpha_1	& \displaystyle\beta_1 & \displaystyle 0 &  \displaystyle 0\\
		\displaystyle 0	& 0 &\displaystyle \alpha_2 &  \displaystyle\beta_2\\
		\displaystyle\frac{1}{2}	& \displaystyle \frac{a}{2 \sqrt{b+p}} &  \displaystyle\frac{-\chi(l_2-l_1,p)}{2} & \displaystyle \frac{-a\chi(l_2-l_1,p)}{2 \sqrt{b+p}} \\
		\displaystyle\frac{-\chi(l_2-l_1,p)}{2}	& \displaystyle\frac{a\chi(l_2-l_1,p)}{2 \sqrt{b+p}} & \displaystyle\frac{1}{2} & \displaystyle \frac{-a}{2 \sqrt{b+p}}
	\end{array}
	\right] \\
	\times \left[
	\begin{array}{cccc}
		\displaystyle G_1(p)\\
		\displaystyle G_2(p)\\
		\displaystyle R(l_1,p) \\
		\displaystyle R(l_2,p)
	\end{array}
	\right];
\end{array}
$$
where $ \chi(x,p)=\exp\left(\frac{-x\sqrt{b+p}}{a} \right) $, and $R(x,p)$ is the Laplace transform of the function
\begin{equation}
	\begin{array}{ll}
		r (x,t)= \displaystyle\frac{\exp(-bt)}{2 a\sqrt{\pi t}}
		\int_{l_1}^{l_2}\varphi(\xi)\left[\exp\left( -\frac{(\xi-x)^2}{4 a^2 t}\right) \right] d\xi\\
		+\displaystyle\frac{1}{2 a\sqrt{\pi }}\int_{0}^{t}d\theta\int_{l_1}^{l_2}\frac{\exp(-b(t-\theta))}{\sqrt{(t-\theta)}}\exp\left( -\frac{(\xi-x)^2}{4 a^2 (t-\theta)}\right)f(\xi,\theta)d\xi.
	\end{array}
	\label{eq18}
\end{equation}

We now return to the problem outlined in equations (\ref{eq06})-(\ref{eq09}). In this context, the source term and the related boundary functions, specifically $g_1(t)$ and $g_2(t)$, are null. Consequently, it is sufficient to assume that the function $ \varphi(x) $ at least satisfies the Dirichlet conditions for  $ x $ in the interval  $[l_1, l_2]$. The unique and exact operational solution of the system (\ref{eq06})-(\ref{eq09}) for the null reaction term ($b=0$), can be reported as:
\begin{equation}
	\begin{array}{ll}
		\displaystyle	U(x,p)= \frac {1}{2\, \sqrt{p}} \,{\exp\left( {{\frac {\left( -l_{{2}}+x \right)  \sqrt{p}}{a}}}\right)} \left[ { U(l_2,p)} \sqrt{p}+a{ U_x(l_2,p)}\right]\\
		\displaystyle+ {\frac {1}{2\, \sqrt{p}} {\exp\left( {-{\frac { \left( x-l_{{1}} \right)  \sqrt{p}}{a}}}\right) }\left[{ U(l_1,p)} \sqrt{p}-a{U_x(l_1,p)}\right]}+R \left( x,p \right),
		\label{eq19}
	\end{array}
\end{equation}
As the source term is null in equation (\ref{eq06}), $R(x,p)$ reduces to the Laplace transform of
\begin{equation}
	r (x,t)= \displaystyle\frac{1}{2 a\sqrt{\pi t}}
	\int_{l_1}^{l_2}\varphi(\xi)\left[\exp\left( -\frac{(\xi-x)^2}{4 a^2 t}\right) \right] d\xi,
	\label{eq20}
\end{equation}
namely,
\begin{equation}
	R(x,p)=\mathcal{L}\{r(x,t)\}=\int_{l_{{1}}}^{l_{{2}}}\! {\frac {\varphi  \left( \xi
			\right) }{2\,a \sqrt{p}}}\, \exp\left( {-{\frac { \left| \xi
				-x \right|  \sqrt{p}}{a}}}\right)  \,{\rm d}\xi,
	\label{eq21}
\end{equation}
with $ \mathcal{L} $ being the Laplace transform operator.  In the present case, it can be written:
\begin{equation}
	\begin{array}{ll}
		\left[
		\begin{array}{cccc}
			\displaystyle U(l_1,p)\\
			\displaystyle U_x(l_1,p)\\
			\displaystyle U(l_2,p) \\
			\displaystyle U_x(l_2,p)
		\end{array}
		\right]	=\left[
		\begin{array}{cccc}
			\displaystyle \alpha_1	& \displaystyle 2a^2 & \displaystyle 0 &  \displaystyle 0\\
			\displaystyle 0	& 0 &\displaystyle \alpha_2 &  \displaystyle 2a^2\\
			\displaystyle\frac{1}{2}	& \displaystyle \frac{a}{2 \sqrt{p}} &  \displaystyle\frac{-\chi(l_2-l_1,p)}{2} & \displaystyle \frac{-a\chi(l_2-l_1,p)}{2 \sqrt{p}} \\
			\displaystyle\frac{-\chi(l_2-l_1,p)}{2}	& \displaystyle\frac{a\chi(l_2-l_1,p)}{2 \sqrt{p}} & \displaystyle\frac{1}{2} & \displaystyle \frac{-a}{2 \sqrt{p}}
		\end{array}
		\right] \\
		\times \left[
		\begin{array}{cccc}
			\displaystyle 0\\
			\displaystyle 0\\
			\displaystyle R(l_1,p) \\
			\displaystyle R(l_2,p)
		\end{array}
		\right].
	\end{array}
	\label{eq22}
\end{equation}
The solutions of system (\ref{eq22}) can be explicitly calculated, and it can be verified that:
\begin{equation}
	\displaystyle\frac{U_x(l_1,p)}{U(l_1,p)}=-{\frac {{ \alpha}_{{1}}}{2 {a}^{2}}},
	\label{eq23}
\end{equation}
and
\begin{equation}
	\displaystyle\frac{U_x(l_2,p)}{U(l_1,p)}=-{\frac {{ \alpha}_{{2}}}{2 {a}^{2}}},
	\label{eq24}
\end{equation}
showing that the boundary conditions of the Laplace analog for the time domain problem (\ref{eq06})-(\ref{eq09}) are satisfied, since $G_1(p)=G_2(p)=0,$ and equations (\ref{eq15})-(\ref{eq16}) reduce to:
\begin{equation}
	\alpha_1 U(l_1,p)+2 {a}^{2} U_x(l_1,p)=0,
	\label{eq25}
\end{equation}
\begin{equation}
	\alpha_2 U(l_2,p)+2 {a}^{2} U_x(l_2,p)=0.
	\label{eq26}
\end{equation}

Therefore, by using the inverse Laplace operator $\mathcal{L}^{-1}$, an exact solution of the problem (\ref{eq01})-(\ref{eq04}) can be derived implicitly from equation (\ref{eq05}) as:
\begin{equation}
	\displaystyle w(x,t)=-2 {a}^{2}\frac{\mathcal{L}^{-1}\{U_x(x,p)\}}{\mathcal{L}^{-1}\{U(x,p)\}},
	\label{eq27}
\end{equation}
where $U(x,p)$ is given by the formula (\ref{eq19}) together with the Laplace domain expression (\ref{eq21}) of $R(x,p) $, and the four solutions of the system (\ref{eq22}). Also, the derivative:
\begin{equation}
	\begin{array}{ll}
		\displaystyle	U_x(x,p)= \frac {1}{2\, a} \,{\exp\left( {{\frac {\left( -l_{{2}}+x \right)  \sqrt{p}}{a}}}\right)} \left[ { U(l_2,p)} \sqrt{p}+a{ U_x(l_2,p)}\right]\\
		\displaystyle- {\frac {1}{2\,a} {\exp\left( {-{\frac { \left( x-l_{{1}} \right)  \sqrt{p}}{a}}}\right) }\left[{ U(l_1,p)} \sqrt{p}-a{U_x(l_1,p)}\right]}+R_x \left(x,p \right),
		\label{eq28}
	\end{array}
\end{equation}
depends on the derivative $R_x(x,p)$ of $R(x,p)$. 
In accordance with the general hypotheses that lead to the unified solution (\ref{eq17}), the formulae (\ref{eq19}) and (\ref{eq28}) are capable of handling smooth as well as generalized solutions that present some discontinuities in their derivatives, even in the case of the space variable (see \cite{An-approximations}). Consequently, relation (\ref{eq27}) is identical to the exact closed-form solution in the time domain, should it exist. In comparison to exact solutions that contain infinite series expressions, solution (\ref{eq27}) represents a ratio of inverse transforms of closed-form expressions. Numerical algorithms for inverse Laplace transforms have been developed that serve as powerful tools for extending the applicability of the Laplace transform technique. An algorithm that has been demonstrated to be particularly efficient is described by de Hoog et al. \cite{Ho-Kn-improved} and has been implemented in Matlab as the function invlap.m by Hollenbeck in 1998. This allows for the graphical illustration of the exact solution (\ref{eq27}) for further comparison.

\section{Illustrations of the Results}\label{sec:3}
The objective of this section is to illustrate the results through the use of specific examples. The curves obtained by the method developed above are compared and discussed with those of exact solutions on two classical examples. The numerical results are obtained using the finite difference method implemented in Matlab via the function pdepe(). 

\subsection{Example 1}\label{subsec:3.1}
The first problem consists to solve equation (\ref{eq01}) with initial condition in the
form:
\begin{equation}
	w_0(x)=w(x,0)=2\,{\frac {{a}^{2}\pi \,\sin \left( \pi \,x \right) }{\sigma+\cos
			\left( \pi \,x \right) }},\quad 0 \le x \le 1,
	\label{eq29}
\end{equation}
where $ \sigma $ is a real parameter. The initial condition of the related problem (\ref{eq06})-(\ref{eq09}) is computed as:
$$
\displaystyle \varphi(x)=u \left( 0,0 \right) {\exp\left( {-\,{\frac {1}{2{a}^{2}}}\int_0^{x}
		\!w_{{0}} \left( y \right) \,{\rm d}y}\right) }={\frac {u \left( 0,0 \right)  \left({\sigma}+\cos \left( {x}
		\,\pi  \right)  \right) }{1+{\sigma}}},
$$
for convenient values of $ \sigma $, among which $ \sigma=2 $ will be taken for graphical illustrations. The following boundary conditions are to be retained:
\begin{equation}
	w\left( 0,t \right)=w\left(1,t \right)=0.
	\label{eq30}
\end{equation}
The exact solution to this problem is represented by a closed-form expression, as detailed in Wood \cite{Wo-exact}, 2014:
\begin{equation}
	w\left( x,t \right)=2\,{\frac {{a}^{2}\pi \,{{\rm e}^{-{\pi }^{2}{a}^{2}t}}\sin \left(
			\pi \,x \right) }{\sigma+{{\rm e}^{-{\pi }^{2}{a}^{2}t}}\cos \left(
			\pi \,x \right) }}.
	\label{eq31}
\end{equation}
\begin{figure}[t]
	\centering
	\includegraphics[width=119mm]{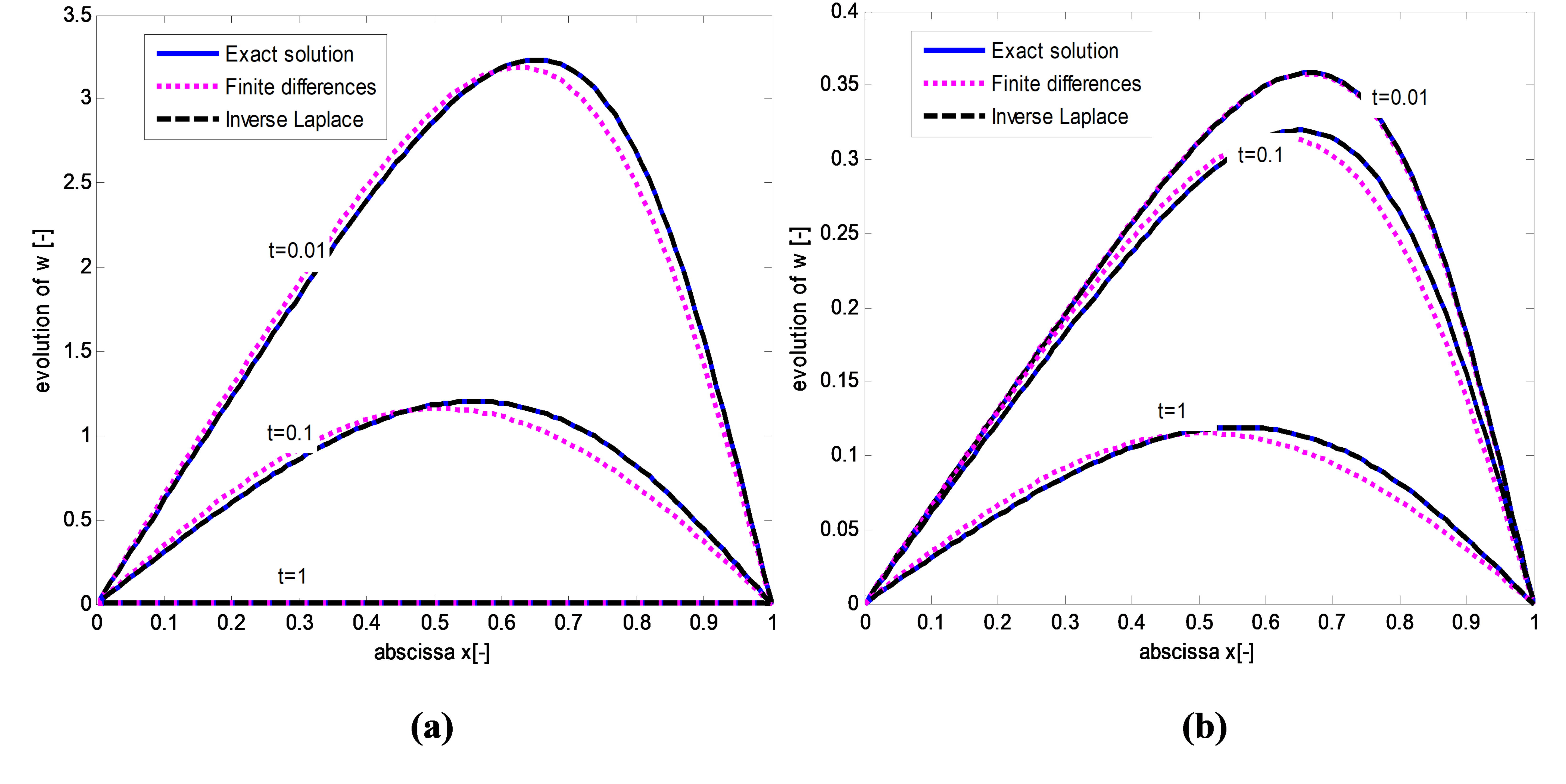}\hfill
	\caption{Solutions $ w(x,t)$  at different times in Example \ref{subsec:3.1}: (\textbf{a}) $ a^2=1,\ \sigma=2 $, (\textbf{b}) $ a^2=0.1,\ \sigma=2$ }
	\label{fig1}
\end{figure}
By assigning the corresponding values to the parameters of system (\ref{eq01})-(\ref{eq04}), that is, $ l_1 =0,\  l_2=1, \ \alpha_1=\alpha_2 =0,$ the Laplace domain function $ R(x,p) $ in formula (\ref{eq21}) can be explicitly calculated. The solutions of system (\ref{eq22}) are as follows:
$$
\left[
\begin{array}{cccc}
	\displaystyle U(l_1,p)\\
	\displaystyle U_x(l_1,p)\\
	\displaystyle U(l_2,p) \\
	\displaystyle U_x(l_2,p)
\end{array}
\right]	= \left[ \begin {array}{c} \displaystyle {\frac { \left( {\pi }^{2}{a}^{2}\sigma
		+p\sigma+p \right) u \left( 0,0 \right) }{3 \left( {\pi }^{2}{a}^{2}+p
		\right) p \left( \sigma+1 \right) }}\\ \noalign{\medskip}0
\\ \noalign{\medskip} \displaystyle
{\frac { \left( {\pi }^{2}{a}^{2}\sigma+p
		\sigma-p \right) u \left( 0,0 \right) }{3 \left( {\pi }^{2}{a}^{2}+p
		\right) p \left( \sigma+1 \right) }}\\ \noalign{\medskip}0
\end {array} \right].
$$
The unique solution of the related problem (\ref{eq06})-(\ref{eq09}) is expressed in the Laplace domain as (see formula (\ref{eq19})):
$$
\displaystyle U(x,p)={\frac {u \left( 0,0 \right)  \left( {\pi }^{2}{a}^{2}\sigma+\cos
		\left( \pi \,x \right) p+p\sigma \right) }{3 \left( \sigma+1 \right) p
		\left( {\pi }^{2}{a}^{2}+p \right) }}.
$$
Its derivative with respect to $ x $ is deduced as:
$$
\displaystyle U_x(x,p)=-{\frac {u \left( 0,0 \right) \sin \left( \pi \,x \right) \pi }{3
		\left( \sigma+1 \right)  \left( {\pi }^{2}{a}^{2}+p \right) }}.
$$
As a result, the exact solution in the time domain of the initial problem under consideration is obtained through the application of the inverse Laplace operator $ \mathcal{L}^{-1} $, as illustrated in formula (\ref{eq27}):
$$
\displaystyle w(x,t)=-2 {a}^{2}\frac{\mathcal{L}^{-1}\{U_x(x,p)\}}{\mathcal{L}^{-1}\{U(x,p)\}}=2\,{\frac {{a}^{2}\pi \,{{\rm e}^{-{\pi }^{2}{a}^{2}t}}\sin \left(
		\pi \,x \right) }{\sigma+{{\rm e}^{-{\pi }^{2}{a}^{2}t}}\cos \left(
		\pi \,x \right) }}.
$$
Figure~\ref{fig1} illustrates the evolution curves of the solution $ w(x,t) $ at different time values for Example \ref{subsec:3.1}.
The curves were obtained by employing a mesh of space step $ \Delta x=0.01$ and time step $ \Delta t=0.001 $. The curves denoted by \lq Finite differences\rq \thinspace represent numerical solutions computed by the Matlab function pdepe(), while the curves designated by \lq Inverse Laplace\rq \thinspace are obtained through numerical calculation of the formula (\ref{eq27}) using the Matlab function invlap.m. The curves align closely with one another in nearly all instances between the exact and numerical inverse Laplace solutions. However, with the exception of the case where $ a=1 $ and $ t=1 $, wherein all the curves collapse to zero, a nonnegative quadratic error ($ L^2 $-norm) can be evaluated between the exact and finite differences solutions.

\subsection{Example 2}\label{subsec:3.2}
Equation (\ref{eq01}) is at $0<x<1 $ for $0<t\le 1 $, and subject to initial condition:
\begin{equation}
	w_0(x)=w(x,0)=\sin(\pi x),
	\label{eq32}
\end{equation}
and to the same homogeneous boundary conditions as in (\ref{eq30}). In this case, the related function $ \varphi $ is calculated as:
$$
\displaystyle \varphi(x)=u \left( 0,0 \right) {{\exp}\left( {{\frac {-1+\cos \left( \pi \,x
				\right) }{2{a}^{2}\pi }}}\right).}
$$
The exact solution contains infinite series in its expression and is denoted by \lq Theoretical solution\rq \thinspace for graphical purposes. It can be written as in \cite{Co-quasilinear}:
\begin{equation}
	w(x,t)= \displaystyle 2 \pi a^2 \frac{\displaystyle \sum_{n=1}^{\infty} c_n   \exp(-n^2 \pi^2 a^2 t) n \sin(n\pi x)}{c_0 +\displaystyle \sum_{n=1}^{\infty} c_n \exp(-n^2 \pi^2 a^2 t) \cos(n\pi x)},
	\label{eq33}
\end{equation}
where the Fourier coefficients are respectively:
$$
c_0=\displaystyle \int_{0}^{1}{\exp}{{\left( \frac {\cos \left( \pi \,x \right) -1}{2 {a}^{2}\pi }\right)}} dx,
$$
and
$$
c_n=\displaystyle 2\int_{0}^{1}{\exp}{{\left( \frac {\cos \left( \pi \,x \right) -1}{2 {a}^{2}\pi }\right)}}\cos(n\pi x) dx,\quad n=1,\,2,\,...
$$
Concerning the $p$-domain expressions, the function $ R(x, p) $ can be written as:
$$
\begin{array}{ll}
	\displaystyle R(x, p)={\frac {u \left( 0,0 \right) }{2a \sqrt{p}}\int_{0}^{1}\!{\exp
			{\left(- {\frac {2\, \left| \xi-x \right|  a\pi\sqrt{p} -\cos \left( \xi
						\,\pi  \right) +1}{2{a}^{2}\pi }}\right) }}\,{\rm d}\xi}\\
	\displaystyle =\frac {u \left( 0,0 \right) }{2a \sqrt{p}} I (0,1,\eta(\left| \xi-x \right|,p),\theta(\xi))\\
	\displaystyle =\frac {u \left( 0,0 \right) }{2a \sqrt{p}} [\eta(x,p)I (0,x,\eta(\xi,p) ,\theta(\xi))+\eta(-x,p)I(x,1,\eta(-\xi,p),\theta(\xi))],
\end{array}
$$
where
$ \displaystyle \eta(\xi,p)=\frac{-\xi\sqrt{p}}{a},$
$\displaystyle \theta(\xi)={{\frac {-1+\cos \left( \pi \,\xi
			\right) }{2{a}^{2}\pi }} },$ and
$\displaystyle I (l_1,l_2,\eta(\xi,p) ,\theta(\xi))=\int_{l_1}^{l_2}\!{\exp
	{\left[{ \eta(\xi,p) \theta(\xi) }\right] }}\,{\rm d}\xi.$
Based on the above notation, the operational solution of system (\ref{eq06})-(\ref{eq09}) can be deduced as:
\begin{equation}
	\begin{array}{ll}
		U(x,p)= K(p)I (0,1,\eta(\xi,p) ,\theta(\xi))[\eta(2-x,p)+\eta(x,p)]
		\\\\
		+K(p)I (0,1,\eta(-\xi,p) ,\theta(\xi))[\eta(2+x,p)+\eta(2-x,p)]+R(x,p),
		\label{eq34}
	\end{array}
\end{equation}
where
$$K(p)=\displaystyle -\frac {u \left( 0,0 \right) }{ 2 a \sqrt{p}\left( {\eta(2,p)}-1 \right)}.
$$
Since $R(x,p)$ is differentiable with respect to $x$, the derivative $U_x(x,p)$ is obtained as:
\begin{equation}
	\begin{array}{ll}
		\displaystyle U_x(x,p)=\frac{\sqrt{p} K(p)}{a}I (0,1,\eta(\xi,p) ,\theta(\xi))[\eta(2-x,p)-\eta(x,p)]\\\\
		\displaystyle+\frac{\sqrt{p} K(p)}{a}I (0,1,\eta(-\xi,p) ,\theta(\xi))[-\eta(2+x,p)+\eta(2-x,p)]+R_x(x,p),
		\label{eq35}
	\end{array}
\end{equation}
where
$$
\begin{array}{ll}
	\displaystyle R_x(x, p) =-\frac {u \left( 0,0 \right) }{2a^2}\\
	\times [\eta(x,p)I (0,x,\eta(-\xi,p) ,\theta(\xi))-\eta(-x,p)I(x,1,\eta(\xi,p),\theta(\xi))].			\end{array}
$$
It can thus be stated that, the exact implicit solution to the problem outlined in Example 2, is given by the equation (\ref{eq27}), with the functions $U(x,p)$ and $U_x(x,p)$ defined in equations (\ref{eq34}) and (\ref{eq35}).
\begin{figure}[t]
	\centering
	\includegraphics[width=119mm]{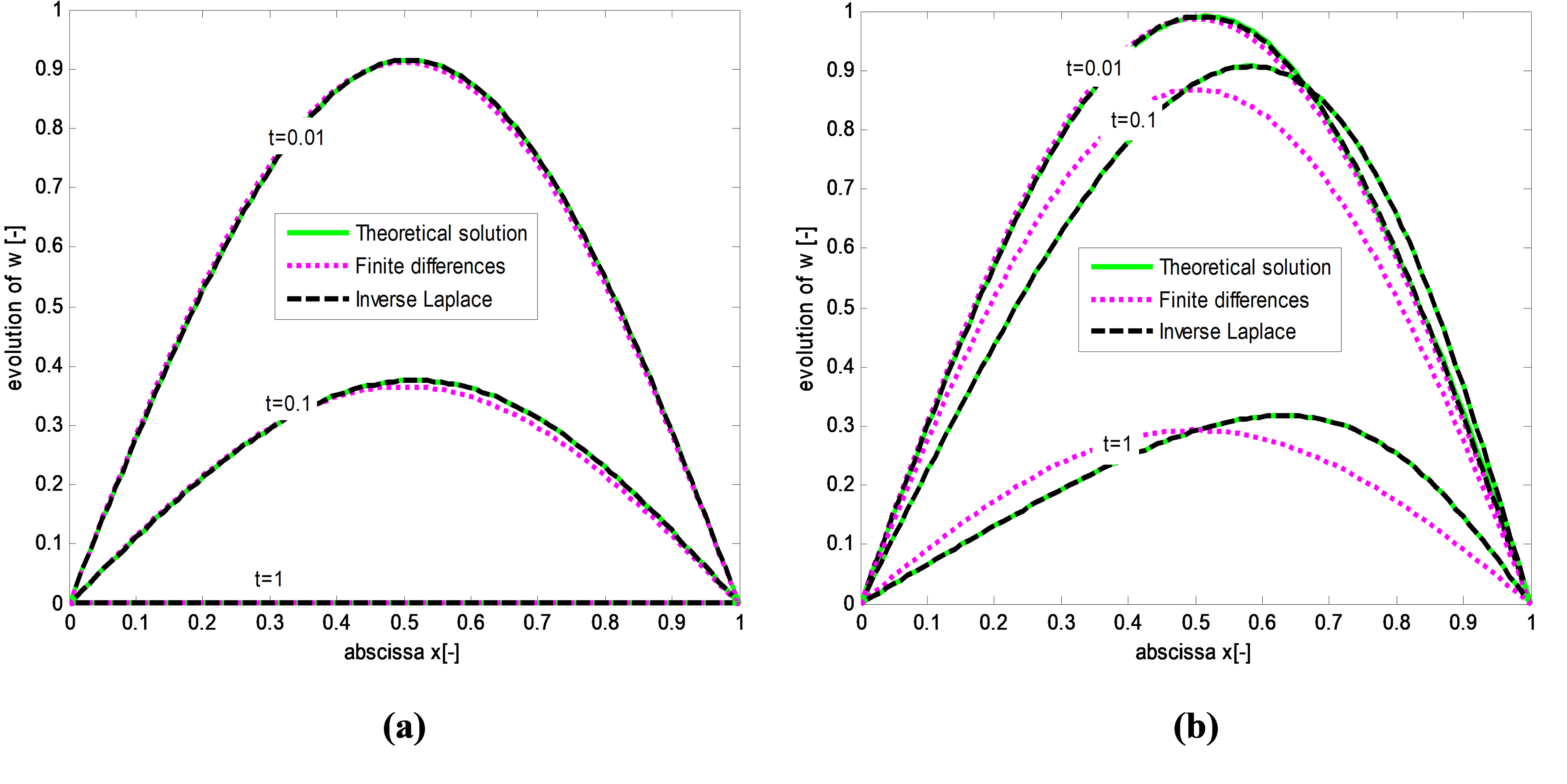}\hfill
	\caption{Solutions $ w(x,t)$  at different times in Example \ref{subsec:3.2}: (\textbf{a}) $ a^2=1$, (\textbf{b}) $ a^2=0.1$}
	\label{fig2}
\end{figure}
The evolution curves of the solution $ w(x,t)$ at different times for the case study outlined in Example \ref{subsec:3.2} are presented in Figure~\ref{fig2}. The mesh cells are identical to those utilized in Example \ref{subsec:3.1}, since reducing the time step does not provide a noticeable improvement. In order to obtain the theoretical solution curve, both infinite series in the ratio (\ref{eq33}) are truncated at their first five terms, as the curves obtained are almost identical when truncating at more terms. The curves obtained for the numerical inverse Laplace transform coincide with those of the theoretical solution for different values of the coefficient $a$. Nevertheless, the script execution time of the numerical inverse Laplace procedure is much reduced compared to that of the theoretical solution (more than 100 times for the first execution). Regarding the finite differences solution curve, except at the very beginning of the process ($t=0.01$), the discrepancy between it and the theoretical solution curve becomes important as the coefficient $a$ decreases, as shown in Figures~\ref{fig2}(a), where $a^2=1$, and \ref{fig2}(b), where $a^2=0.1$. In brief, the numerical inverse Laplace procedure, by using the exact implicit formula (\ref{eq27}), is notable for its computational efficiency.


\section{Conclusion}
This study has facilitated the precise calculation of an exact solution for the viscous Burgers' equation on bounded intervals by combining the use of the Hopf-Cole transformation with the Laplace integral transform. In instances where an explicit closed-form analytic inverse in the time domain is unavailable, the exact implicit solution has been demonstrated to be more efficient than classical numerical and semi-analytical methods. This is due to the availability of highly efficient algorithms for numerical inverse Laplace transforms. The Burgers' equation has a variety of applications in numerous fields of physics and applied mathematics. The equation is a valuable tool for understanding a wide range of nonlinear and diffusive phenomena across various scientific disciplines, offering insights at both the theoretical and practical levels. The closed-form exact solution in the Laplace or time domains will be employed as a point of reference for numerical and semi-analytical methods. In future work, we will investigate the application of a comparable solution methodology to the forced Burgers' equation in one and more spatial dimensions. Furthermore, the present solution methodology will be examined in relation to other approaches for solving nonlinear evolution equations, including Painlevé analysis and Hirota's bilinear technique.


\begin{thebibliography}{00}

\bibitem{Ba-nonlinear}
Ba\c{s}han, A.: Nonlinear dynamics of the Burgers' equation and numerical experiments. Math. Sci., Springer. 16(2), 183--205 (2022). https://doi.org/10.1007/s40096-021-00410-8

\bibitem{Ku-Ku-evolutionary}
Kumar S, Kumar A, Mohan B.: Evolutionary dynamics of solitary wave profiles and abundant analytical solutions to a (3+1)-dimensional burgers system in ocean physics and hydrodynamics. J Ocean Eng Sci. 8(1), 1--14 (2023). https://doi.org/10.1016/j.joes.2021.11.002

\bibitem{Kh-Iq-mathematical}
Khan, M., Iqbal, Z., Ahmed, A.: A mathematical model to examine the heat transport features in Burgers fluid flow due to stretching cylinder. Journal of Thermal Analysis and Calorimetry, 147(1), 827--841 (2022).  https://doi.org/10.1007/s10973-020-10224-w

\bibitem{Ha-Gu-anomalies1}
Hafez, M.M., Guo, W.H.: Some anomalies of numerical simulation of shock waves. Part I: inviscid flows. Comp. Fluids, 28(4-5), 701--719 (1999). https://doi.org/10.1016/S0045-7930(98)00051-6

\bibitem{Ha-Gu-anomalies2}
Hafez, M.M., Guo, W.H.: Some anomalies of numerical simulation of shock waves. Part II: effect of artificial and real viscosity. Comp. Fluids, 28(4-5), 721--739 (1999). https://doi.org/10.1016/S0045-7930(98)00052-8

\bibitem{Zh-Wu-continuous}
Zhai, C., Wu, W.: A continuous traffic flow model considering predictive headway variation and preceding vehicle's taillight effect. Phys. A Stat. Mech. Appl., 584, 126364. (2021). https://doi.org/10.1016/j.physa.2021.126364

\bibitem{Zh-al-cooperative}
Zhai, C., Wu, W., Xiao, Y.: Cooperative car-following control with electronic throttle and perceived headway errors on gyroidal roads. Appl. Math. Modelling, 108, 770--786 (2022).
https://doi.org/10.1016/j.apm.2022.04.010

\bibitem{Gu-Ji-modeling}
Gu, J., Jing, Y.: Modeling of wave propagation for medical ultrasound: a review. IEEE Trans. Ultrason. Ferroelectr. Freq. Control., 62(11), 1979--1992 (2015). 10.1109/TUFFC.2015.007034

\bibitem{Ch-al-sonification}
Chiroiu, V., Munteanu, L., Ioan, R., Dragne, C., Majercsik, L.: Using the sonification for hardly detectable details in medical images. Scientific Reports, 9(1), 17711 (2019).

\bibitem{Ru-He-quadratically}
Rudenko, O.V., Hedberg, C.M.: The quadratically cubic Burgers equation: an exactly solvable nonlinear model for shocks, pulses and periodic waves. Nonlinear Dyn., 85, 767--776 (2016). https://doi.org/10.1007/s11071-016-2721-5

\bibitem{Ku-al-dynamics}
Kumar, R., Kumar, M., Tiwari, A.K.: Dynamics of some more invariant solutions of (3+ 1)-Burgers system. Int. J. Comput. Meth. Eng. Sci. Mech., 22(3), 225--234 (2021). https://doi.org/10.1080/15502287.2021.1916693

\bibitem{Bo-Aw-systematic}
Bonkile, M.P., Awasthi, A., Lakshmi, C., Mukundan, V., Aswin, V.S.: A systematic literature review of Burgers' equation with recent advances. Pramana - J. Phys. 90, 1--21 (2018)

\bibitem{Ho-partial}
Hopf, E.: The partial differential equation $u_t+u u_x = \mu u_{xx}$. Comm. Pure Appl. Math. 3, 201--230 (1950). http://dx.doi.org/10.1002/cpa.3160030302

\bibitem{Co-quasilinear}
Cole, J.D.: On a quasilinear parabolic equations occurring in aerodynamics. Quart. Appl. Math. 9(3), 225--236 (1951). https://www.jstor.org/stable/43633894

\bibitem{Be-Pl-table}
Benton, E.R., Platzman, G.W.: A table of solutions of the one-dimensional Burgers' equations. Quart. Appl. Math. 30(2), 195--212 (1972). https://doi.org/10.1090/qam/306736


\bibitem{Dh-Ka-contemporary}
Dhawan, S., Kapoor, S., Kumar, S., Rawat, S.: Contemporary review of techniques for the solution of nonlinear Burgers equation. J. Comput. Sci., 3(5), 405--419. (2012). https://doi.org/10.1016/j.jocs.2012.06.003


\bibitem{Na-Ma-mathematical}
Naghipour, A., Manafian, J.: Application of the Laplace Adomian decomposition and implicit methods for solving Burgers' equation. TWMS J. Pure Appl. Math. 6(1), 68--77 (2015)


\bibitem{Ze-Ch-mathematical}
Zeidan, D., Chau, C.K., Lu, T.T., Zheng, W.Q.: Mathematical studies of the solution of Burgers' equations by Adomian decomposition method. Math. Methods Appl. Sci. 43(5), 2171--2188 (2020). https://doi.org/10.1002/mma.5982

\bibitem{La-Ya-approximate}
Lal, D., Yadav, M.: Approximate analytical solution of one dimensional nonlinear Burger's equation using Homotopy Perturbation method. J. Algebr. Stat. 13(3), 5462--5469 (2022)

\bibitem{Os-al-approximate}
Osman, M.S., Baleanu, D., Adem, A.R., Hosseini, K., Mirzazadeh, M., Eslami, M.: Double-wave solutions and Lie symmetry analysis to the (2+ 1)-dimensional coupled Burgers equations. Chin. J. Phys. 63, 122--129 (2020). https://doi.org/10.1016/j.cjph.2019.11.005


\bibitem{Gu-al-fifth}
Guo, Y., Shi, Y.F.,  Li, Y.M.: A fifth-order finite volume weighted compact scheme for solving one-dimensional Burgers’ equation. Appl. Math. Comput. 281, 172--185 (2016). https://doi.org/10.1016/j.amc.2016.01.061


\bibitem{Se-accurate}
Seydao$\mathrm{\check{g}}$lu, M.: An accurate approximation algorithm for Burgers’ equation in the presence of small viscosity. J. Comput. Appl. Math. 344, 473-481 (2018). https://doi.org/10.1016/j.cam.2018.05.063


\bibitem{Bi-Am-exact}
Biazar, J., Aminikhah, H.: Exact and numerical solutions for non-linear Burger's equation by VIM. Math. Comput. Modelling 49(7), 1394--1400 (2009). https://doi.org/10.1016/j.mcm.2008.12.006


\bibitem{An-approximations} Anani, K.: Analytical approximations in short times of exact operational solutions to reaction-diffusion problems on bounded intervals. Appl. Appl. Math. 19(1), 1--27 (2023) 


\bibitem{Her-partial} Herron, I.H., Foster, M.R.: Laplace transform methods. In: Herron I.H., Foster, M.R. (eds.) Partial Differential Equations in Fluid Dynamics, pp. 148-182. Cambridge University Press, Cambridge (2008)

\bibitem{De-integral} Debnath, L., Bhatta, D.: Laplace transforms and their basic properties. In: Debnath, L., Bhatta, D. (eds.) Integral Transforms and Their Applications, pp. 143--196. CRC press, Boca Raton (2014)


\bibitem{Ho-Kn-improved}
de Hoog, F.R., Knight J.H., Stokes, A.N.: An improved method for numerical inversion of Laplace transforms. SIAM J. Sci. Stat. Comput. 3(3), 357--366 (1982). https://doi.org/10.1137/0903022


\bibitem{Wo-exact}
Wood, W.L.: An exact solution for Burgers' equation. Commun. Numer. Methods Eng. 22(7), 797--798 (2006). https://doi.org/10.1002/cnm.850




\end{thebibliography}
\end{document}